\documentclass[reqno]{amsart}
\usepackage{amssymb,hyperref}
\usepackage[dvips]{graphicx}

\theoremstyle{plain}
\newtheorem{thm}{Theorem}
\newtheorem{lem}{Lemma}
\newtheorem{cor}{Corollary}

\theoremstyle{definition}
\newtheorem{rem}{Remark}

\renewcommand{\Re}{\mathrm{Re}}
\renewcommand{\Im}{\mathrm{Im}}

\title
[Extensions of Nunokawa lemma for argument properties]
{Extensions of Nunokawa lemma \\
for argument properties}

\author{Hitoshi Shiraishi}
\address{Hitoshi Shiraishi \newline
Department of Mathematics \newline
Kinki University \newline
Higashi-Osaka, Osaka 577-8502, Japan}
\email{step\_625@hotmail.com}


\date{}

\begin{document}

\begin{abstract}
Let $\mathcal{H}[a_0,n]$ be the class of functions $p(z)=a_0+a_nz^n+\cdots$ which are analytic in the open unit disk $\mathbb{U}$.
For $p(z)\in\mathcal{H}[1,2]$, M. Nunokawa, S. Owa, N. Uyanik and H. Shiraishi (Math. Comput. Modelling. {\bf 55} (2012), 1245--1250) have shown some theorems for argument properties.
The object of the present paper is to discuss some extensions of Nnnokawa lemma and its applications for argument properties.
\end{abstract}

\maketitle

\section{Introduction}

\

Let $\mathcal{H}[a_0,n]$ denote the class of functions $p(z)$ of the form
$$
p(z)
= a_0 +\sum^{\infty}_{k=n} a_k z^k
$$
which are analytic in the open unit disk $\mathbb{U}=\{z \in \mathbb{C}:|z|<1\}$ for some $a_0\in\mathbb{C}$ and a positive integer $n$.
If $p(z)\in\mathcal{H}[a_0,n]$ satisfies
$$
|\arg(p(z))|
< \frac{\pi}{2}\mu
\qquad(z\in\mathbb{U})
$$
for some real $0<\mu\leqq1$, then we say that $p(z)$ belongs to the class $\mathcal{STH}[a_0,n](\mu)$.

\

The basic tool in proving our results is the following lemma due to S. S. Miller and P. T. Mocanu \cite{m1ref2} (also \cite{d7ref2}).

\

\begin{lem} \label{jack} \quad
Let the function $w(z)$ definded by
$$
w(z)
= a_n z^n + a_{n+1}z^{n+1} + a_{n+2}z^{n+2} + \cdots
\qquad(n=1,2,3,\cdots)
$$
be analytic in $\mathbb{U}$ with $w(0)=0$.
If $\left|w(z)\right|$ attains its maximum value on the circle $|z|=r$ at a point $z_{0}\in\mathbb{U}$,
then there exists a real number $m \geqq n$ such that
$$
\frac{z_{0}w'(z_{0})}{w(z_{0})}
= m.
$$
\end{lem}

\

\section{Main results}

\

Applying Lemma \ref{jack},
we derive the following result.

\

\begin{thm} \label{p02thm1} \quad
Let $p(z)\in\mathcal{H}[a_0,n]$ for some real $a_0>0$ and suppose that there exists a point $z_0\in\mathbb{U}$ such that
$$
\Re(p(z))
> 0
\quad for \quad |z|<|z_0|
$$
and $p(z_0) = \beta i$ is a pure imaginary number for some real $\beta\neq0$.

Then we have
$$
\frac{z_0 p'(z_0)}{p(z_0)}
= il
$$
where
$$
l
\geqq \frac{n}{2}\left( \frac{a_0}{\beta}+\frac{\beta}{a_0} \right)
\geqq n
$$
if $\beta>0$ and
$$
l
\leqq \frac{n}{2}\left( \frac{a_0}{\beta}+\frac{\beta}{a_0} \right)
\leqq -n
$$
if $\beta<0$.
\end{thm}

\

\begin{proof} \quad
Let us put
$$
w(z)
= \frac{a_0-p(z)}{a_0+p(z)}
= c_nz^n + c_{n+1}z^{n+1} + c_{n+2}z^{n+2} + \cdots
\qquad(z\in\mathbb{U}).
$$

Then, we have that $w(z)$ is analytic in $|z|<|z_0|$, $w(0)=0$, $|w(z)|<1$ for $|z|<|z_0|$ and
$$
|w(z_0)|
= \left| \frac{a_0^2-\beta^2-2a_0\beta i}{a_0^2+\beta^2} \right|
= 1.
$$

From Lemma \ref{jack},
we obtain
$$
\frac{z_0 w'(z_0)}{w(z_0)}
= \frac{-2 a_0 z_0 p'(z_0)}{a_0^2 - \{ p(z_0) \}^2}
= \frac{-2 a_0 z_0 p'(z_0)}{a_0^2 + \beta^2}
= m
\qquad(m \geqq n).
$$

This shows that
$$
z_0 p'(z_0)
= -\frac{m}{2}\left( a_0 + \frac{\beta^2}{a_0} \right)
\qquad(m \geqq n).
$$

From the fact that $z_0 p'(z_0)$ is a real number and $p(z_0)$ is a pure imaginary number,
we can put
$$
\frac{z_0p'(z_0)}{p(z_0)}
= il
$$
where $l$ is a real number.

For the case $\beta>0$,
we have
\begin{align*}
l
&= \Im\left( \frac{z_0p'(z_0)}{p(z_0)} \right) \\
&= \Im\left( -z_0p'(z_0)\frac{1}{\beta}i \right) \\
&= \frac{m}{2}\left( a_0 + \frac{\beta^2}{a_0} \right)\frac{1}{\beta} \\
&\geqq \frac{n}{2}\left( a_0 + \frac{\beta^2}{a_0} \right)\frac{1}{\beta} \\
&= \frac{n}{2}\left( \frac{a_0}{\beta}+\frac{\beta}{a_0} \right)
\geqq n
\end{align*}
and for the case $\beta<0$,
we get
\begin{align*}
l
&= \Im\left( \frac{z_0p'(z_0)}{p(z_0)} \right) \\
&= \Im\left( -z_0p'(z_0)\frac{1}{\beta}i \right) \\
&= \frac{m}{2}\left( a_0 + \frac{\beta^2}{a_0} \right)\frac{1}{\beta} \\
&\leqq \frac{n}{2}\left( a_0 + \frac{\beta^2}{a_0} \right)\frac{1}{\beta} \\
&= \frac{n}{2}\left( \frac{a_0}{\beta}+\frac{\beta}{a_0} \right)
\leqq -n.
\end{align*}

This completes our proof.
\end{proof}

\

Putting $a_0=1$ in Theorem \ref{p02thm1},
we have Corollary \ref{p02cor1}.

\

\begin{cor} \label{p02cor1} \quad
Let $p(z)\in\mathcal{H}[1,n]$ and suppose that there exists a point $z_0\in\mathbb{U}$ such that
$$
\Re(p(z))
> 0
\quad for \quad |z|<|z_0|,
$$
$\Re(p(z_0))=0$ and $p(z_0)\neq0$.

Then we have
$$
\frac{z_0 p'(z_0)}{p(z_0)}
= il
$$
where $l$ is a real and $|l|\geqq n$.
\end{cor}

\

From Theorem \ref{p02thm1},
we get Theorem \ref{p02thm2}.

\

\begin{thm} \label{p02thm2} \quad
Let $p(z)\in\mathcal{H}[a_0,n]$ for some real $a_0<0$ and suppose that there exists a point $z_0\in\mathbb{U}$ such that
$$
\Re(p(z))
< 0
\quad for \quad |z|<|z_0|
$$
and $p(z_0) = \beta i$ is a pure imagenary number for some real $\beta\neq0$.

Then we have
$$
\frac{z_0 p'(z_0)}{p(z_0)}
= il
$$
where
$$
l
\geqq \frac{n}{2}\left( \frac{a_0}{\beta}+\frac{\beta}{a_0} \right)
\geqq n
$$
if $\beta<0$ and
$$
l
\leqq \frac{n}{2}\left( \frac{a_0}{\beta}+\frac{\beta}{a_0} \right)
\leqq -n
$$
if $\beta>0$.
\end{thm}

\

\begin{proof} \quad
We put the function
$$
q(z)
= -p(z)
\qquad(z\in\mathbb{U}),
$$
$q(z)$ satisfies the assumption of Theorem \ref{p02thm1} and using Theorem \ref{p02thm1},
we get the result of Theorem \ref{p02thm2}.
\end{proof}

\

\section{Applications of Theorem \ref{p02thm1}}

\

Using Theorem \ref{p02thm1},
we obtain following result.

\

\begin{thm} \label{p02thm3} \quad
If $p(z)\in\mathcal{H}[a_0,n]$ for some real $a_0>0$ and a positive integer $n\geqq2$ satisfies $p(z)\neq0$ for $z\in\mathbb{U}$ and
$$
|\arg(p(z)-zp'(z))|
< \arctan(n\mu) - \frac{\pi}{2}\mu
\qquad(z\in\mathbb{U})
$$
for some real number $0<\mu\leqq\sqrt[]{\dfrac{2n-\pi}{n^2\pi}}<1$,
then $p(z)\in\mathcal{STH}[a_0,n](\mu)$.
\end{thm}

\

\begin{proof} \quad
Let us consider
$$
q(z)
= (p(z))^\frac{1}{\mu}
= a_0^\frac{1}{\mu}+c_nz^n+c_{n+1}z^{n+1}+\cdots
\qquad(z\in\mathbb{U})
$$
suppose that $p(z)$ satisfies
$$
|\arg (p(z))|
< \frac{\pi}{2}\mu
\qquad(|z|<|z_0|)
$$
and
$$
|\arg(p(z_0))|
= \frac{\pi}{2}\mu
\qquad(z_0\in\mathbb{U}).
$$

Then, the function $q(z)$ satisfies
$$
\Re(q(z))
>0
\qquad(|z|<|z_0|)
$$
and $\Re(q(z_0))=0$ with $q(z_0)\neq0$.

Applying Theorem \ref{p02thm1},
we have
\begin{equation}
\frac{z_0q'(z_0)}{q(z_0)}
= \frac{1}{\mu}\frac{z_0p'(z_0)}{p(z_0)}
= il
\label{p02thm3eq1}
\end{equation}
where
\begin{equation}
l
\geqq \frac{n}{2}\left( \frac{\beta}{a_0^\frac{1}{\mu}}+\frac{a_0^\frac{1}{\mu}}{\beta} \right)
\geqq n
\qquad(q(z_0)=i\beta)
\label{p02thm3eq2}
\end{equation}
and
\begin{equation}
l
\leqq -\frac{n}{2}\left( \frac{\beta}{a_0^\frac{1}{\mu}}+\frac{a_0^\frac{1}{\mu}}{\beta} \right)
\leqq -n
\qquad(q(z_0)=-i\beta)
\label{p02thm3eq3}
\end{equation}
for $\beta>0$ which satisfies
$$
q(z_0)
= (p(z_0))^\frac{1}{\mu}
= \pm i\beta.
$$

For such $p(z)$,
if $\arg(p(z_0))=\dfrac{\pi}{2}\mu$,
we have that
\begin{align*}
\arg(p(z_0)-z_0p'(z_0))
&= \arg\left( p(z_0)\left( 1-\frac{z_0p'(z_0)}{p(z_0)} \right) \right) \\
&= \frac{\pi}{2}\mu + \arg(1-il\mu) \\
&\leqq \frac{\pi}{2}\mu + \arg(1-in\mu) \\
&= \frac{\pi}{2}\mu - \arctan(n\mu) \\
&= -\left( \arctan(n\mu)-\frac{\pi}{2}\mu \right)
\end{align*}
which contradicts the condition in the theorem.

If $\arg(p(z_0))=-\dfrac{\pi}{2}\mu$,
we get
\begin{align*}
\arg(p(z_0)-z_0p'(z_0))
&= \arg(1-il\mu) - \frac{\pi}{2}\mu \\
&\geqq \arg(1+in\mu) - \frac{\pi}{2}\mu \\
&= \arctan(n\mu)-\frac{\pi}{2}\mu
\end{align*}
which contradicts the condition in the theorem.

This implies that there is no $z_0\in\mathbb{U}$ such that
$$
|\arg(p(z))|
< \frac{\pi}{2}\mu
\qquad(|z|<|z_0|)
$$
and
$$
|\arg(p(z_0))|
= \frac{\pi}{2}\mu.
$$

Thus $p(z)$ satisfies
$$
|\arg(p(z))|
< \frac{\pi}{2}\mu
$$
for all $z\in\mathbb{U}$, that is, $p(z)\in\mathcal{STH}[a_0,n](\mu)$.
\end{proof}

\

\begin{rem} \label{p02rem1} \quad
Consider the function
$$
g(x)
= \arctan(nx)-\frac{\pi}{2}x
$$
for some real number $0<x\leqq\sqrt[]{\dfrac{2n-\pi}{n^2\pi}}\leqq\sqrt[]{\dfrac{4-\pi}{4\pi}}$ and some positive integer $n\geqq2$,
we get
$$
g'(x)
= \frac{2n-\pi-n^2\pi x^2}{2(1+n^2x^2)}
\geqq 0
\qquad\left( 0<x\leqq\sqrt[]{\dfrac{2n-\pi}{n^2\pi}} \right).
$$

From $g(0)=0$,
we know that $g(z)$ is a simple increasing function and $g(x)>0$.
\end{rem}

\

When $n=2$ and $a_0=1$ in Theorem \ref{p02thm3},
we have Corollary \ref{p02cor2} due to M. Nunokawa, S. Owa, N. Uyanik and H. Shiraishi \cite{ds3ref0}.

\

\begin{cor} \label{p02cor2} \quad
If $p(z)\in\mathcal{H}[1,2]$ satisfies $p(z)\neq0$ for $z\in\mathbb{U}$ and
$$
|\arg(p(z)-zp'(z))|
< \arctan(2\mu) - \frac{\pi}{2}\mu
\qquad(z\in\mathbb{U})
$$
for some real number $0<\mu\leqq\sqrt[]{\dfrac{4-\pi}{4\pi}}$,
then $p(z)\in\mathcal{STH}[1,2](\mu)$.
\end{cor}

\

Also, using Theorem \ref{p02thm1},
we obtain the following theorem.

\

\begin{thm} \label{p02thm4} \quad
If $p(z)\in\mathcal{H}[a_0,n]$ for some real $a_0>0$ and a positive integer $n\geqq1$ satisfies $p(z)\neq0$ for $z\in\mathbb{U}$ and
\begin{equation}
\left| \arg \left( p(z)-\frac{zp'(z)}{p(z)} \right) \right|
< \arctan \left( \dfrac{\rho(\mu) \cos \dfrac{\pi\mu}{2} }{1-\rho(\mu) \sin \dfrac{\pi\mu}{2} } \right) - \frac{\pi}{2}\mu
\qquad(z\in\mathbb{U})
\label{p02thm4eq1}
\end{equation}
for some real number $\mu$ $(0<\mu<\mu_0)$,
where
\begin{equation}
\rho(\mu)
= \frac{n\mu}{2 a_0} \left( \left( \frac{1+\mu}{1-\mu} \right)^\frac{1-\mu}{2} + \left( \frac{1-\mu}{1+\mu} \right)^\frac{1+\mu}{2} \right)
\label{p02thm4eq2}
\end{equation}
and some real number $0<\mu_0<1$ satisfies
$$
\rho(\mu_0)
= \sin \frac{\pi\mu_0}{2},
$$
then $p(z)\in\mathcal{STH}[a_0,n](\mu)$.
\end{thm}

\

\begin{proof} \quad
The condition (\ref{p02thm4eq1}) implies, in particular, that $p(z)\neq0$ for $z\in\mathbb{U}$.
We consider that there exists a point $z_0\in\mathbb{U}$ such that
$$
|\arg(p(z))|
< \frac{\pi}{2}\mu
\qquad(|z|<|z_0|)
$$
and
$$
|\arg(p(z_0))|
= \frac{\pi}{2}\mu.
$$

By using the same prosess of Theorem \ref{p02thm3},
we obtain the equation (\ref{p02thm3eq1}) and we can write
$$
\frac{z_0p'(z_0)}{p(z_0)}
= i \mu l,
$$
where $l$ is a real number which satisfies the equation (\ref{p02thm3eq2}) and (\ref{p02thm3eq3}).

For the case $\arg(p(z_0))=-\dfrac{\pi}{2}\mu$,
applying the same method as the proof of Theorem \ref{p02thm3},
we have that
\begin{align*}
&\arg\left( p(z_0)-\frac{z_0p'(z_0)}{p(z_0)} \right) \\
&= \arg(p(z_0)) + \arg \left( 1-\frac{z_0p'(z_0)}{p(z_0)}\frac{1}{p(z_0)} \right) \\
&= -\frac{\pi}{2}\mu + \arg \left( 1-\frac{i \mu l}{(-i\beta)^\mu} \right) \\
&= -\frac{\pi}{2}\mu + \arg \left( 1-e^{i\frac{1+\mu}{2}\pi}\frac{\mu}{\beta^\mu}l \right) \\
&\geqq -\frac{\pi}{2}\mu + \arg \left( 1+e^{i\frac{1+\mu}{2}\pi}\frac{n\mu}{2\beta^\mu} \left( \frac{\beta}{a_0^\frac{1}{\mu}}+\frac{a_0^\frac{1}{\mu}}{\beta} \right) \right) \\
&= -\frac{\pi}{2}\mu + \arg \left( 1+e^{i\frac{1+\mu}{2}\pi}\frac{n\mu}{2a_0} \left( \left( \frac{\beta}{a_0^\frac{1}{\mu}} \right)^{1-\mu} + \left( \frac{\beta}{a_0^\frac{1}{\mu}} \right)^{-1-\mu} \right) \right).
\end{align*}

On the other hand, let us put
$$
g(x)
= x^{1-\mu} + x^{-1-\mu}
\qquad \left( x=\frac{\beta}{a_0^\frac{1}{\mu}}>0 \right).
$$

Then, by easy calculation, we have
$$
g'(x)
= (1-\mu)x^{-\mu} + (-1-\mu)x^{-2-\mu}
\qquad(x>0)
$$
and $g(x)$ takes the minimum value at
$$
x
= \sqrt[]{\frac{1+\mu}{1-\mu}}.
$$

Therefore, we have
\begin{align*}
&\arg\left( p(z_0)-\frac{z_0p'(z_0)}{p(z_0)} \right) \\
&\geqq -\frac{\pi}{2}\mu + \arg \left( 1+e^{i\frac{1+\mu}{2}\pi}\frac{n\mu}{2a_0} \left( \left( \frac{\beta}{a_0^\frac{1}{\mu}} \right)^{1-\mu} + \left( \frac{\beta}{a_0^\frac{1}{\mu}} \right)^{-1-\mu} \right) \right) \\
&\geqq -\frac{\pi}{2}\mu + \arg \left( 1+e^{i\frac{1+\mu}{2}\pi}\frac{n\mu}{2a_0} \left( \left( \frac{1+\mu}{1-\mu} \right)^\frac{1-\mu}{2} + \left( \frac{1-\mu}{1+\mu} \right)^\frac{1+\mu}{2} \right) \right) \\
&= \arctan \left( \dfrac{\rho(\mu) \sin \left( \dfrac{1+\mu}{2}\pi \right) }{1+\rho(\mu) \cos \left( \dfrac{1+\mu}{2}\pi \right) } \right) - \frac{\pi}{2}\mu \\
&= \arctan \left( \dfrac{\rho(\mu) \cos \dfrac{\pi\mu}{2} }{1-\rho(\mu) \sin \dfrac{\pi\mu}{2} } \right) - \frac{\pi}{2}\mu
\end{align*}
for the equation (\ref{p02thm4eq2}),
which contradicts the condition (\ref{p02thm4eq1}).

If $\arg(p(z_0))=\dfrac{\pi}{2}\mu$,
then we see that
\begin{align*}
&\arg\left( p(z_0)-\frac{z_0p'(z_0)}{p(z_0)} \right) \\
&= \frac{\pi}{2}\mu + \arg \left( 1-\frac{i \mu l}{(i\beta)^\mu} \right) \\
&\leqq \frac{\pi}{2}\mu - \arctan \left( \dfrac{\rho(\mu) \sin \left( \dfrac{1+\mu}{2}\pi \right) }{1+\rho(\mu) \cos \left( \dfrac{1+\mu}{2}\pi \right) } \right) \\
&= \frac{\pi}{2}\mu - \arctan \left( \dfrac{\rho(\mu) \cos \dfrac{\pi\mu}{2} }{1-\rho(\mu) \sin \dfrac{\pi\mu}{2} } \right)
\end{align*}
for the equation (\ref{p02thm4eq2}),
which also contradicts the condition (\ref{p02thm4eq1}).

This shows that there is no $z_0\in\mathbb{U}$ such that
$$
|\arg(p(z))|
< \frac{\pi}{2}\mu
\qquad(|z|<|z_0|)
$$
and
$$
|\arg(p(z_0))|
= \frac{\pi}{2}\mu.
$$

Therefore, $p(z_0)$ satisfies $|\arg(p(z))|<\dfrac{\pi}{2}\mu$ for all $z\in\mathbb{U}$,
completing the proof.
\end{proof}

\

\begin{rem} \label{p02rem2} \quad
Let the function
$$
g(x)
= \arctan \left( \dfrac{\rho(x) \cos \dfrac{\pi x}{2} }{1-\rho(x) \sin \dfrac{\pi x}{2} } \right) - \frac{\pi}{2} x
$$
where $\rho(x)$ satisfies (\ref{p02thm4eq2}) for some positive integer $n\geqq2$ and some real $a_0>0$.

We know $g(x)>0$ for all $x$ $(0<x\leqq\mu_0)$ where some real number $0<\mu_0<1$ satisfies
$$
\rho(\mu_0)
= \sin \frac{\pi\mu_0}{2}.
$$

Because,
we have the following figure for the function $g(x)$.

\begin{figure*}[h]
\begin{center}
\includegraphics[clip,scale=0.4]{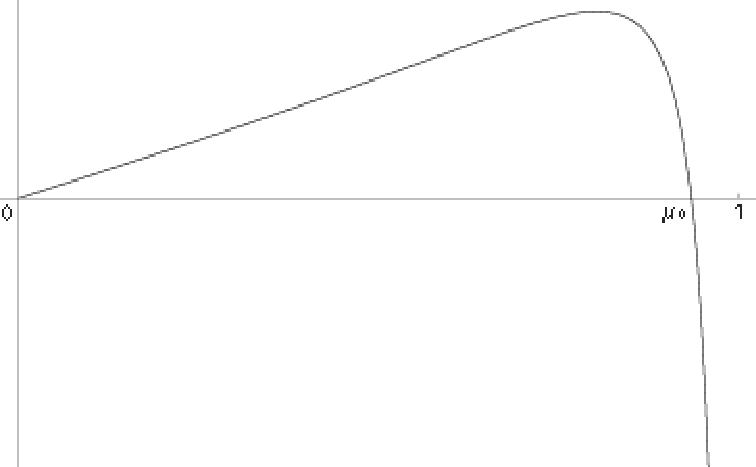}
\end{center}
\caption{The image of $g(x)$.}
\end{figure*}

\end{rem}

\

Considering $n=2$ and $a_0=1$ in Theorem \ref{p02thm4},
we have the following corollary due to M. Nunokawa, S. Owa, N. Uyanik and H. Shiraishi \cite{ds3ref0}.

\

\begin{cor} \label{p02cor3} \quad
If $p(z)\in\mathcal{H}[1,2]$ satisfies $p(z)\neq0$ for $z\in\mathbb{U}$ and
\begin{align*}
\left| \arg \left( p(z)-\frac{zp'(z)}{p(z)} \right) \right|
&< \arctan \left( \dfrac{\rho(\mu) \cos \dfrac{\pi\mu}{2} }{1-\rho(\mu) \sin \dfrac{\pi\mu}{2} } \right) - \frac{\pi}{2}\mu \\
&= \arctan \left( \dfrac{\rho(\mu) \sin \left( \dfrac{1+\mu}{2}\pi \right) }{1+\rho(\mu) \cos \left( \dfrac{1+\mu}{2}\pi \right) } \right) - \frac{\pi}{2}\mu
\qquad(z\in\mathbb{U})
\end{align*}
for some real number $\mu$ $(0<\mu<\mu_0)$,
where
$$
\rho(\mu)
= \mu \left( \left( \frac{1+\mu}{1-\mu} \right)^\frac{1-\mu}{2} + \left( \frac{1-\mu}{1+\mu} \right)^\frac{1+\mu}{2} \right)
$$
and some real number $0<\mu_0<1$ satisfies
$$
\rho(\mu_0)
= \sin \frac{\pi\mu_0}{2},
$$
then $p(z)\in\mathcal{STH}[1,2](\mu)$.
\end{cor}

\

\end{document}